\documentclass[12pt]{article}
\usepackage{amssymb,amsfonts,amsmath,amsthm}
\usepackage{epsfig}
\newtheorem{thm}{Theorem}[section]

\newtheorem{lem}{Lemma}[section]

\newtheorem{defn}[thm]{Definition}

\makeatletter \@addtoreset{equation}{section}

\def\pf{\noindent {\it Proof.\ }}
\def\qed{\hfill \rule{4pt}{7pt}}

\begin{document}
\begin{center}
{\large\bf Derangements and Relative Derangements of Type $B$}\\
[7pt]
William Y.C. Chen$^1$ and Jessica C.Y. Zhang$^2$\\
[5pt]

Center for Combinatorics, LPMC-TJKLC\\ Nankai University, Tianjin
300071, P. R. China\\ [5pt]

$^1$chen@nankai.edu.cn, $^2$zhangcy@mail.nankai.edu.cn

\end{center}
\vskip 3mm

\begin{abstract}
By introducing the notion of relative derangements of type $B$, also
called signed relative derangements, which are defined in terms of
signed permutations, we obtain a type $B$ analogue of the well-known
relation between relative derangements and the classical
derangements. While this fact can be proved by using the principle
of inclusion and exclusion, we present a combinatorial
interpretation with the aid of the intermediate structure of signed
skew derangements.
\end{abstract}



\section{Introduction}

  A
 derangement on a set $[n]=\{1, 2,\cdots, n\}$ is a permutation
 $\pi=\pi_{1}\pi_{2}\cdots
 \pi_{n}$ such that $\pi_{i}\neq i$ for all $i\in [n]$. A
 relative derangement $\pi_{1}\pi_{2}\cdots \pi_{n}$ on [n] is a
 permutation such that $\pi_{i+1}\neq \pi_{i}+1$ for $1\leq i \leq
 n-1$. Let $Q_n$ denote the number of relative derangements on
 $[n]$, and let $D_{n}$ denote the number of the derangements on $[n]$.
The following relation is well-known, see Brualdi \cite[Theorem
6.5.1]{Brualdi}, or Andreescu and Feng \cite[Example 6.11]{And}:
\begin{equation} \label{me-a}
Q_n=D_{n}+D_{n-1}.
\end{equation}

 A combinatorial interpretation of (\ref{me-a}) has been
obtained by Chen \cite{Chen} based on the intermediate structure of
skew derangements, which are equivalent to the generalized
derangements as studied by Hanson, Seyffarth and Weston \cite{HSW}
and Wang \cite{Wang}. The main objective of this paper is to present
a type $B$ analogue of (\ref{me-a}). This goal is achieved by
introducing the notion of signed relative derangements, or relative
derangements of type $B$. The concept of derangements of type $B$ is
introduced by Chow \cite{Chow}. A {\it signed permutation}  $\pi$ on
$[n]$  can be viewed as a bijection on the set $\{\bar{1}, \cdots,
\bar{n}, 1, \cdots, n\}$ such that $\pi(\bar{i})=\overline{\pi(i)}$.
Intuitively, a signed permutation on $[n]$ is just an ordinary
permutation $\pi_1\pi_2\cdots \pi_n$ with some elements associated
with a bar $-$. For example, $3\, \bar{2}\, \bar{5}\, 1\, \bar{4}$
is a signed permutation on $\{1, 2,3,4,5\}$. The set of signed
permutations on $[n]$ is often denoted by $B_n$. The following order
relation is often imposed on the elements of signed permutations for
$B_n$, see, for example, Shareshian and Wachs \cite{Wachs}:
\begin{equation}\label{order}
\bar{1}<\bar{2}<\cdots<\bar{n}<1<2<\cdots<n.
\end{equation}

According to the above ordering, for the above signed permutation
$3\, \bar{2}\, \bar{5}\, 1\, \bar{4}$, $3$ is the largest element
and $\bar{2}$ is the smallest. We recall the following definition of
derangements of type $B$.

\begin{defn}
A derangement of type $B$ on $[n]$ is a signed permutation
$\pi_{1}\pi_{2}\cdots\pi_{n}$ such that $\pi_{i}\neq i$, for all
$i\in [n]$.
\end{defn}

For example,  $3\, \bar{2}\, \bar{5}\, 1\, \bar{4}$ is a derangement
in $B_{5}$, whereas $3\, 2\, \bar{4}\, 1\, \bar{5}$  has a fixed
point $2$. Let $D_{n}^{B}$ denote the number of derangements of type
$B$ on $[n]$. It is not hard to derive the following formula by
using the principle of inclusion-exclusion [4, Chapter 2]:

\begin{eqnarray}
 D_{n}^{B}=n!\sum_{k=0}^{n}\frac{(-1)^{k}\cdot 2^{n-k}}{k!}
\end{eqnarray}
In fact, it is also a consequence of the $q$-analogue given by Chow
\cite{Chow}.

We now give the definition of relative derangements of type $B$ on
$[n]$, or signed relative derangements, for short.

\begin{defn} A
relative derangement of type $B$ on $[n]$ is a signed permutation on
$[n]$ such that $i$ is not followed by $i+1$, and $\bar{i}$ is not
followed by $\overline{i+1}$, for $1 \leq i \leq n-1$.
\end{defn}

For example,  $3\, 2\, \bar{4}\, 1\, \bar{5}$ is a relative
derangement in $B_{5}$, while $4\, 1\, 5\, \bar{2}\, \bar{3}$ is
not. Let $Q_n^{B}$ be the number of relative derangements of type
$B$. Our main result is the following type $B$ analogue of the above
relation (\ref{me-a}).

\begin{thm}\label{1} For $n\geq 2$, we have
\begin{eqnarray} \label{qd}
Q_n^{B}=D_{n}^{B}+D_{n-1}^{B}.
\end{eqnarray}
\end{thm}

The first few values of $Q_n^B$ starting with $Q_1^B$ are given
below:
\[ 2,\ 6,\ 34,\ 262,\ 2562,\ \cdots\]
In accordance with the relation (\ref{qd}), we adopt the convention
that $D_0^B=1$.

 One way to prove the above result for $Q_n^B$  and $D_n^B$ is to
derive the following formula for $Q_n^B$ by using the principle of
inclusion-exclusion:
\begin{eqnarray} \label{qnb}
 Q_n^{B}=n!\cdot2^{n}+\sum_{k=1}^{n-1}(-1)^{k}\cdot{n-1\choose
 k}\cdot(n-k)!\cdot2^{n-k}.
\end{eqnarray}
However, the details of the algebraic proof will be omitted.
Instead, we will provide a combinatorial proof by introducing the
structure of signed skew derangements.

\section{Signed Skew Derangements}

In this section, we first introduce the notion of signed skew
derangements and establish a correspondence between signed relative
derangements and signed skew derangements. Then we give a
characterization of signed permutations that correspond to signed
skew derangements. Then we show how to transform a signed skew
derangement into a signed derangements. This leads to a
combinatorial interpretation of the relation (\ref{qd}).

Recall that a skew derangement $f$ on $[n]$  is a bijection from
$[n]$ onto $\{0, 1,\cdots, n-1\}$ with $f(i)\neq i$ for any $i\in
[n]$, see \cite{Chen}. For signed permutations, we will define
signed skew derangements, or  skew derangements of type $B$. Let us
begin with the definition of a signed set on $[n]$. A signed set on
$[n]$ can be considered the underlying set of a signed permutation.
In other words, a signed set on $[n]$ is just the set $[n]$ with
some elements bearing bars.
 For example,
$X=\{1, \bar{2}, 3, 4, \bar{5}\}$ is a signed set on $\{1, 2, 3,
4,5\}$.

Given a signed set $X$ on $[n]$, we denote by $X-1$ the signed set
obtained from $X$ by subtracting $1$ from each element in $X$, where
we define the subtraction for barred elements by the rule
\begin{equation}\label{rule1}
 \overline{i} -  1 = \overline{i-1}.\end{equation}
Conversely, the addition to a  barred element is given by
\begin{equation}\label{rule2}
 \overline{i} + 1 = \overline{i+1}.
 \end{equation}

\begin{defn}
Let $X$ be a signed set on $[n]$.  A signed skew derangement on
$[n]$ is a bijection $f$ from $X$ to $Y=X-1$ such that $f(x)\neq x$
for any $x\in X$, where $x$ may be a barred element.
\end{defn}

For example, let $n=2$, $X=\{\bar{1}, 2\}$ and $Y=\{\bar{0}, 1\}$.
Then there are two signed skew derangements from $X$ to $Y$:
$f_{1}(\bar{1})=\bar{0}$, $f_{1}(2)=1$ and $f_{2}(\bar{1})=1$,
$f_{2}(2)=\bar{0}$. The following theorem establishes a bijection
between signed relative derangements and signed skew derangements.

\begin{thm}\label{4}
There is a one-to-one correspondence between the set of signed
relative derangements on $[n]$ and the set of signed skew
derangements on $[n]$.
\end{thm}

\pf First, given a signed relative derangement
$\pi=\pi_{1}\pi_{2}\cdots \pi_{n}$ on $[n]$, we proceed to construct
a signed skew derangement $f$ on $[n]$. Let $u$ be the maximum
element in the signed permutation $\pi_{1}\pi_{2}\cdots \pi_{n}$
with respect to the order (\ref{order}). Note that in the case of
signed permutations, the maximum element is not necessarily the
element $n$.
 Suppose that $\pi_{k}=u$. Let us consider
the segment $\pi_{1}\pi_{2}\cdots \pi_{k}$. Define
$$f(\pi_{1})=\pi_{2}-1,\ \ f(\pi_{2})=\pi_{3}-1,\ \ \cdots, f(\pi_{k-1})=\pi_{k}-1,\ \
f(\pi_{k})=\pi_{1}-1,\ \ $$ subject to the above subtraction rule
(\ref{rule1}) if an element $\pi_t$ is a barred element.

By the definition of signed relative derangement, we claim that $f$
satisfies the condition of a signed skew derangement with respect to
the elements $\pi_1, \pi_2, \ldots, \pi_k$, namely,
\[ f(\pi_1)\neq \pi_1, \quad f(\pi_2)\neq \pi_2, \quad \ldots, \quad
f(\pi_k)\neq \pi_k.\]

For any $r=1, 2, \cdots,k-1$, since $\pi$ is a signed relative
derangement, in view of the addition operation (\ref{rule2}) we see
that $\pi_{r+1}\neq \pi_r+1$ no matter whether $\pi_r$ is a barred
element or not. So we have
$$f(\pi_{r})=\pi_{r+1}-1\neq \pi_{r}$$
for $r=1, 2, \cdots,k-1$. We now consider $\pi_k$.  Since $\pi_{k}$
is the maximum element of $\pi$, we find $\pi_{1}-1\neq \pi_{k}$.
This implies that $f(\pi_{k})=\pi_{1}-1\neq \pi_{k}$.

Now we can repeat the above procedure for the remaining sequence
$\sigma=\pi_{k+1} \pi_{k+2} \cdots \pi_{n}$. The next step is still
to choose the maximum element $\pi_{t}$ in $\sigma$, then assign the
images of $f$ for the elements $\pi_{k+1}, \pi_{k+2},\ldots,
\pi_{t}$. If there are still elements left, we may iterate this
procedure until $f$ is completely determined.

It remains to construct the inverse procedure.
 Given a signed skew
derangement $f$ on $[n]$, we aim to find
 the corresponding
signed relative derangement.

Suppose $f$ is a bijection from a signed set $X$ to $X-1$. The first
step is to determine $\pi_1$. Assume that $u$ is the maximum element
in $X$ with respect to the order (\ref{order}). Then we set
$\pi_1=f(u)+1$, subject to the above addition rule (\ref{rule2}) if
$f(u)$ is a barred element. Suppose $\pi_{r}$ is already located. If
$\pi_{r}\neq u$, then we set $\pi_{r+1}=f(\pi_{r})+1$, using the
above rule (\ref{rule2}) if
 $f(\pi_r)$ is a barred element, and repeat this process until
 we reach a step when  $\pi_{k}=u$ for some $k$.

 At this point, we have obtained the segment $i_{1}i_{2}\cdots i_{k}$.
Since $f(i_r)\neq i_r$, we see that $i_{r+1}\neq i_{r}+1$, for $r=1,
\cdots, k$. If $k<  n$, then  we may choose the maximum element in
the remaining elements in $X$ after removing the elements $i_1, i_2,
\ldots, i_k$, and iterate the above procedure until we obtain the
desired signed relative derangement. Thus, we have shown that our
construction is a bijection. \qed

For example, the signed relative derangement $\bar{7}\, 8\, 6\,
\bar{1}\, \bar{5}\, \bar{3}\, 4\, 2$ corresponds to the following
signed skew derangement:
\begin{align*}
&f(\bar{7})=8-1=7,\quad f(8)=\bar{7}-1=\bar{6},\quad
f(6)=6-1=5,\quad
f(\bar{1})=\bar{5}-1=\bar{4}, \\[10pt]
&f(\bar{5})=\bar{3}-1=\bar{2},\quad  f(\bar{3})=4-1=3,\quad
f(4)=\bar{1}-1=\bar{0},\quad f(2)=2-1=1.
\end{align*}

We now turn our attention to a combinatorial interpretation of the
fact that the number of signed skew derangements on $[n]$ equals
$D_n^B +D_{n-1}^B$. As the first step, we give a characterization of
signed permutations on $\{0, 1,\ldots, n-1\}$ that correspond to
signed skew derangements on $[n]$. Let us consider bijections from a
signed set $X$ on $[n]$ to $X-1$. Assume that the elements of $X$
are arranged by the increasing order of their underlying elements,
say, $X=\{\sigma_1, \sigma_2, \ldots, \sigma_n\}$.  It is easy to
observe the fact that a bijection $f$ from $X$ to $X-1$ is determine
by the signed permutation $\pi=\pi_1\pi_2\cdots\pi_n$, where $\pi_i=
f(\sigma_i)$. In fact, this is a bijection, because for any signed
permutation $\pi$ on $\{0, 1, \ldots, n-1\}$, the elements $\{
\pi_1, \pi_2, \ldots, \pi_n\}$ determines the signed set $X-1$,
which in turn determines $X$. Hence the map $f$ from $X$ to $X-1$ is
easily constructed. The signed permutation $\pi$ is called the {\it
representation} of $f$.

For the above signed skew derangement $f$, we have \[ X=\{\sigma_1,
\sigma_2, \ldots, \sigma_8\} =\{\bar{1}, 2, \bar{3}, 4, \bar{5}, 6,
\bar{7}, 8\}\]
 and $\pi=\bar{4}\, 1\, 3\, \bar{0}\, \bar{2}\, 5\,
7\, \bar{6}$.

The following lemma gives a characterization of signed permutations
which are representations of signed skew derangements. A bar
associated with an element is intuitively considered as a sign.
Moreover, for a signed permutation $\pi=\pi_1\pi_2\cdots \pi_n$, an
element $\pi_i$ is called a {\it fixed point} if $\pi_i=i$, whereas
it is called a {\it signed fixed point} if $\pi_i=i$ or $\bar{i}$.
As will be seen, signed fixed points play an important role in
establishing the correspondence between signed skew derangements and
signed derangements.

\begin{lem} Let $\pi$ be a signed permutation on $\{0, 1, \ldots,
n-1\}$, and let $X$ and $f$ be the signed set and the bijection from
$X$ to $X-1$ determined by $\pi$.  Then $f$ has a fixed point if
$\pi$ has a signed fixed point $\pi_{i}$, and $i-1$ and $i$ have the
same sign in $\pi$.
\end{lem}

The above lemma can be restated as follows. A signed permutation
$\pi$ is a representation of a signed skew derangement if and only
if  $\pi_i=i$ implies that $\overline{i-1}$ appears in $\pi$, and
$\pi=\bar{i}$ implies that $i-1$ appears in $\pi$.

\noindent {\it Proof.} Let $\pi$ be a signed permutation on $\{0, 1,
2,\ldots, n-1\}$. Let $f$ be a bijection from $X$ to $X-1$ such that
 $\pi$ is the representation of $f$. Then $X-1$ is determined by
the entries of $\pi$. Hence $X$ is uniquely determined by $\pi$. Let
$\sigma_1, \sigma_2, \ldots, \sigma_n$ be the elements of $X$
arranged in the increasing order of the underlying elements of $X$.
If $f$ has a fixed point, say, $f(x)=x$, for some $x=\sigma_i$. Then
we have $\sigma_i=i$ or $\bar{i}$, and $f(\sigma_i)=\sigma_i=\pi_i$.
Since $f$ is a bijection from $X$ to $X-1$, $\sigma_i$ is a barred
element if and only if $i-1$ is a barred element. Thus, we conclude
that $\pi_i$ and $i-1$ have the same sign. This completes the proof.
\qed

The above characterization indicates that signed skew derangements
can be viewed as an intermediate structure between signed relative
derangements and signed derangements.  Using this characterization
of representations of signed skew derangements on $[n]$, we first
consider a class of such signed permutations that are in one-to-one
correspondence with signed derangements on $[n-1]$.

\begin{lem} \label{L-S-A}
There is a bijection between the set of representations
 of signed skew
derangements on $[n]$ that are of the form $\pi=\pi_1\pi_2\cdots
\pi_{n-1}0$  and the set of signed derangements on $[n-1]$.
\end{lem}

For example, there are five signed derangements on $\{1,2\}$:
$\bar{1}\bar{2}$, $21$, $2\bar{1}$, $\bar{2}1$, $\bar{2}\bar{1}$. In
the meantime, there are five representations signed skew
derangements on $\{1, 2,3\}$ that are of the form $\pi_1\pi_20$:
$\bar{1}\bar{2}0$, $210$, $2\bar{1}0$, $\bar{2}10$,
$\bar{2}\bar{1}0$. As in this example, special attention should be
paid to the signed derangement $\bar{1}\bar{2}$ with signed fixed
points, and to the representation $\bar{1}\bar{2}0$ which also have
signed fixed points. In general, we can establish a correspondence
as given in the following proof.

\noindent {\it Proof.} Let $\pi=\pi_1\pi_2\cdots \pi_{n-1}0$ be a
representation of a signed skew derangement on $[n]$. We aim to
construct a signed derangement on $[n-1]$ from $\pi$. If
$\pi_1\pi_2\cdots \pi_{n-1}$ has no signed fixed point, then it is
automatically the desired signed derangement.

We now consider that case when there are some signed fixed points,
namely, there exist some $i$ such that $\pi_i=i $ or $\bar{i}$.
Taking the signed fixed point $\pi_i$
 with  minimum index $i$,
we observe that whether $\pi_i$ has a bar or not is determined
solely by the appearance of $i-1$ in the sense that it is a barred
element or an unbarred element. Iterating this argument, we may
deduce that the signed fixed points are uniquely determined by the
remaining elements in $\pi$. Hence we may always put $\bar{i}$ as
the signed fixed points in order to obtain a signed derangement.

Conversely, given a signed derangement $\tau=\tau_1\tau_2\cdots
\tau_{n-1}$, we may identify the signed fixed points $\tau_i$. By
the same argument as in the previous paragraph, we can determine the
signed fixed points according to the characterization of
representations of signed skew derangements so that the resulting
signed permutation on $\{0, 1, \ldots, n-1\}$ corresponds to a
signed skew derangement. This completes the proof.  \qed

For example, consider the signed skew derangement $f$ on
$\{1,2,\ldots, 8\}$ which has the following representation
\[ f(1)\, f(2)\, f(\bar{3})\, f(\bar{4})\, f(5)\,
f(\bar{6})\, f(\bar{7})\, f(8)=\bar{6}\, \bar{2}\, 1\, 4\, \bar{3}\,
7\, \bar{5}\, 0 .\] It corresponds to the signed derangement
$\bar{6}\, \bar{2}\, 1\, \bar{4}\, \bar{3}\, 7\, \bar{5}$ on $\{1,
2, \ldots, 7\}$.

To complete the combinatorial proof of Theorem \ref{1}, it suffices
to consider the second case for the representations of signed skew
derangements. The following lemma is concerned with this case.

\begin{lem} There is a one-to-one correspondence between
representations  $\pi=\pi_1\pi_2\cdots \pi_n$ of signed skew
derangements on $[n]$ with $\pi_n\neq 0$ and signed derangements on
$[n]$.
\end{lem}

For example, there are five representation $\pi=\pi_1\pi_2$ of
signed skew derangements on $\{1,2\}$ with $\pi_2\neq 0$: $01$,
$0\bar{1}$, $\bar{0}1$, $\bar{0}\bar{1}$, $1\bar{0}$.

\noindent {\it Proof.} First, we show that from a representation
$\pi=\pi_1\pi_2\cdots \pi_n$ of a signed skew derangement with
$\pi_n\neq 0$ we can construct a signed derangement
$\tau=\tau_1\tau_2\cdots \tau_n$. If there is no signed fixed point
in $\pi$, then we can replace $0$ or $\bar{0}$ by $n$ or $\bar{n}$
in $\pi$ depending whether $0$ or $\bar{0}$ appears. Since
$\pi_n\neq 0$, we have $\tau_n\neq n$ and so  the resulting signed
permutation is a signed derangement on $[n]$.

Otherwise, there are some signed fixed points $\pi_i$ $(1\leq i \leq
n-1)$, namely, $\pi_i=i$ or $\bar{i}$. Using the same argument as in
the proof of Lemma \ref{L-S-A}, we see that the signed fixed points
are completely determined by the remaining elements in the signed
permutation. So we may set all the signed fixed points to barred
elements in $\pi$.  Finally, we may replace $0$ by $n$ or $\bar{0}$
by $\bar{n}$ to get a signed derangement $\tau$ on $[n]$.

It is easy to see that the above procedure is reversible. This
completes the proof. \qed

For example, consider the signed skew derangement $f$ on $\{1, 2,
\ldots, 8\}$ which has the following representation
\[ f(\bar{1})\, f(2)\, f(\bar{3})\, f(4)\, f(\bar{5})\,
f(6)\, f(\bar{7})\, f(8)=\bar{4}\, 1\, 3\, \bar{0}\, \bar{2}\, 5\,
7\, \bar{6}.\] The corresponding  signed derangement turns out to be
$\bar{4}\, 1\, \bar{3}\, \bar{8}\, \bar{2}\, 5\, \bar{7}\, \bar{6}$.

Combining the preceding two lemmas leads to a combinatorial
interpretation of Theorem \ref{1}. To conclude this paper, we remark
that our bijection between signed relative derangements  and signed
skew derangements can be restricted to ordinary permutations. Hence
the classical relation (\ref{me-a}) is  a consequence of Theorem
\ref{1}.

\vskip 8pt
 \noindent {\bf Acknowledgments.}  This work was supported
by the 973 Project, the PCSIRT Project of the Ministry of Education,
the Ministry of Science and Technology, and the National Science
Foundation of China.


\end{document}